\documentclass{amsart}
\usepackage{amssymb}

\usepackage{tikz-cd}

\def\1{^{-1}}

\def\CP{{\mathbf C\mathbf P}}
\def\MU{{\mathbf M\mathbf U}}
\def\MSU{{\mathbf M\mathbf S \mathbf U}}

\def\MSp{{\mathbf M\mathbf S \mathbf p}}

\newtheorem{theorem}{Theorem}[section]
\newtheorem{lemma}[theorem]{Lemma}

\theoremstyle{definition}

\theoremstyle{remark}

\newcommand{\BU}{\mathop{\mathrm{BU}}}

\numberwithin{equation}{section}

\begin{document}



\commby{}

\title[]{Complex cobordism $\MU^*[1/2]$ modulo $\MSU^*[1/2]$ and related genera}

\begin{abstract}
This paper presents a commutative complex oriented cohomology theory with coefficients the quitient ring  of complex cobordism $\MU^*[1/2]$ modulo 
the ideal generated by any subsequence  of any polynomial generators in special unitary cobordism $\MSU^*[1/2]$ 
viewed as elements in $\MU^*[1/2]$ by forgetful map.    
\end{abstract}

\author{Malkhaz Bakuradze }
\address{Faculty of exact and natural sciences, A. Razmadze Math. Institute, Iv. Javakhishvili Tbilisi State University, Georgia }
\email{malkhaz.bakuradze@tsu.ge } 
\thanks{The author was supported by Shota Rustaveli NSF grant FR-21-4713}

\maketitle

\section{Introduction}

The theory of complex cobordism $\MU^*(-)$ and special unitary cobordism  $\MSU^*(-)$ play an important role in cobordism theory \cite{Stong}.

In particular the coefficient ring $\MSU_*$, localized away from 2, is torsion free
$$\MSU_*[1/2]=\mathbb{Z}[1/2][x_2,x_3,\cdots],\,\,\,\,|x_i|=2i$$
and $SU$-structure forgetful homomorphism is the inclusion 
$$\MSU_*[1/2]\subset \MU_*[1/2]=\mathbb{Z}[1/2][a_1,a_2,a_3,\cdots].$$

With this in mind the generators $x_i$ in $\MSU_*[1/2]$ can be treated as elements in $\MU_*[1/2]$. In particular $x_i$ is a $SU$ manifold iff all Chern numbers of $x_i$ having factor $c_1$ are zero. 
Then we have to check the main Chern number $s_i(x_i)$ for Novikov's criteria  \cite{NOV} for the membership of the set of polynomial generators in  
$\MSU^*[1/2]$. 

\bigskip

In Section 3 we prove

\begin{lemma}
	\label{regular}
 Let  $\Gamma=\{x_n, \, n\geq 2\}$ be a sequence of any polynomial generators in $\MSU_*[1/2]$ viewed as elements in $\MU_*[1/2]$ by forgetful map. Then
 
 i) $\, \Gamma$ is regular. 
 
 ii) any subsequence $\Sigma$ of  $\, \Gamma$ is regular.
\end{lemma}

\bigskip

By using Baas-Sullivan theory of cobordism with singularities we deduce
    
\begin{theorem}
\label{h}
For any subset $\Sigma$ of the set of polynomial generators in $\MSU_*[1/2]$  
there exists a commutative complex oriented cohomology  
$\emph{h}_{\Sigma}^*(-):=\MU^*_{\Sigma}[1/2](-)$,
with the coefficient ring $\MU_*[1/2]/(\Sigma)$.
\end{theorem}

Note that the examples in Lemma \ref{regular} include   $\Sigma=\{x_3, x_4, \cdots \}$. Then Theorem \ref{h} gives abelian cohomology \cite{BUSATO} after localized away from 2. 

Another example was done in \cite{B2,B4}. If $\Sigma=\{x_n,\, n\geq 5 \}$ then $h^*_{\Sigma}$ is a commutative complex oriented cohomology with scalar ring isomorphic to 
the ring of coefficients of the universal Buchstaber formal group law localized away from 2.

In Section 4 we discuss some related genera.

\bigskip

\section{Preliminaries}

Let, as above,  $x_i$ be polynomial generators in $\MSU_*[1/2]$ viewed  as elements in $\MU_*[1/2]$ using the forgetful map. Let 
$$F_U=\sum \alpha_{ij}x^iy^j$$ 
be the universal formal group law of complex cobordism and let $\CP_i$ be the cobordism class of complex projective spae of dimension $2i$.

Then $x_2, x_2$ and $x_4$ are as follows \cite{B5}

\begin{align}
	\label{234}
	&x_2=\CP_2-\frac{9}{8}\CP_1^2,	
	&&x_3=-\alpha_{22},
	&&x_4=-\alpha_{23}-\frac{3}{2}x_3\CP_1.
\end{align}


These are the $SU$ manifolds as all Chern numbers of $x_i$ having factor $c_1$ are zero. Then we have to check the main Chern number $s_n(x_n)$ for Novikov's criteria  \cite{NOV} for the membership of the set of polynomial generators in  
$\MSU^*[1/2]$.

In particular, $SU$ manifold $x_n$ is a polynomial
generator if and only if 
\begin{align}
	\label{Novikov}
	s_n(x_n)=
	\begin{cases} \pm 2^kp &\mbox{if } n=p^l,\,\,\,p\,\,\, \mbox{is odd prime, }\\
		\pm 2^kp &\mbox{if } n+1=p^l,\,\,\, p\,\,\,\mbox{is odd prime, } \\
		\pm 2^k& \mbox{otherwise}.
	\end{cases}
\end{align}

By Euclid's algorithm for the natural numbers $m_1,m_2, \cdots, m_k$  one can find integers $\lambda_1, \lambda_2, \cdots \lambda_k $ such that 
\begin{equation}
	\label{lambda}
	\lambda_1m_1+\lambda_2m_2+\cdots  +\lambda_km_k=gcd(m_1,m_2,\cdots ,m_k).
\end{equation}

Let 
\begin{equation}
	\label{d}
	d(m)=gcd \bigg \{ \binom{m+1}{1},\binom{m+1}{2},\cdots ,\binom{m+1}{m-1}   \big|\,\, m\geq 1 \bigg \}.
\end{equation}

By \cite{KU} one has 
\begin{equation*}
	d(m)=
	\begin{cases} p, &\mbox{if $m+1=p^s$ for some prime  $p$,}\\
		1, &\mbox{otherwise}.
	\end{cases}
\end{equation*}

  The elements

\begin{equation}
	\label{e_n}
	e_{m}=\lambda_1\alpha_{1\,m}+\lambda_2\alpha_{2\,m-1}+\cdots +\lambda_m\alpha_{m\,1}	
\end{equation}  
are multiplicative generators in $\MU_*$ \cite{?}

Let $m\geq 4$ and let $\lambda_2,\cdots ,\lambda_{m-2}$ are such integers that

\begin{equation}
	\label{d_2}
	d_2(m):=\sum_{i=2}^{m-2}\lambda_i\binom{m+1}{i}=gcd \bigg \{ \binom{m+1}{2},\cdots ,\binom{m+1}{m-2}  \bigg \}.
\end{equation}

Then by \cite{BU-U} Lemma 9.7 one has for $m\geq 3$

\begin{align} 
	\label{d_2d_1} 
	&d_2(m)=d(m)d(m-1).	
\end{align}

By definition the coefficient ring of the universal abelian formal group law $F_{Ab}$ is the quotient ring
\begin{equation}
\label{Ab}
\Lambda_{Ab}=\MU_*/I_{Ab}, \text { where } I_{Ab}=(\alpha_{ij}, i,j >1),
\end{equation}
where $\alpha_{ij}$ are the coefficients of the universal formal group law $F_U$.

Let us apply Euclid's algorithm for the Chern numbers $$s_{m-1}(\alpha_{i,m-i})=\binom{m+1}{i}$$
in \eqref{d_2} and let
$$y_{k}=\sum_{i=2}^{k-1}\lambda_i\alpha_{i\, k+1-i},\,\,\,k\geq 3. $$   

Recall from  \cite{B-KH}, \cite{BUSATO} the constriction of abelian cohomology $h^*{Ab}(-)$. By \cite{BUSATO} one has $I_{Ab}=(y_k, \,\,k\geq 3)$. Moreover 
the sequence $Y=\{y_k, \,k\geq 3\}$ is regular. Using Baas-Sullivan theory this defines $h^*{Ab}(-)$, which is commutative as all obstructions vanish.

\section{Proofs}

To prove Lemma \ref{regular} i) let as above $\Sigma$ be the ideal generated by any sequence   $\Sigma=\{x_n,\, n\geq 5 \}$ of polynomial generators in $\MSU_*$ and let $(Y)$ be the ideal generated by the regular sequence $Y=\{y_k, \,k\geq 3\}$ defined above.  It is clear that the latter remains regular after localized away from 2. 

One has $(\Sigma)\subset (Y)$ because the abelian formal group law is a specialization of the Buchstaber formal group law, which has a greater kernel ideal. By \eqref{234} we have  $(\Gamma)\subset (x_2,Y)$. 

It is clear that  $x_2$ is a polynomial generator in $\MU_*[1/2]$ therefore the  ideal $(x_2,Y)$ is also regular.

It suffices to prove the reverse inclusion 
$((x_2,Y))\subset \Gamma)$. It can proved by induction to follow the proof of Proposition 5.3 in \cite{B5} and taking into account the structure of the quotient ring $\MU_*/(Y)$ \cite{B-KH}. We leave the detailed proof to the reader.


ii) For the regularity of any subsequence $\{x_{i_1}, \cdots x_{i_n}\}$ of $\Gamma$ we need three facts: if $R$ is a graded ring and the members of a regular sequence are homogeneous of positive degree, then any permutation of a regular sequence is  regular \cite {BRUNS}.  Then it is clear by definition that  any first $n$ members of a regular sequence form a regular sequence and vice versa (including infinite case).  

\qed

The proof of Theorem \ref{h} uses 
the Sullivan-Baas construction \cite{BA}  of cobordism with singularities. The regular sequence  
$\Sigma$ gives a cohomology theory $\MU^*_{\Sigma}(-)$ which by regularity 
of the ideal $(\Sigma)$ has a scalar ring 
$$\MU^*_{\Sigma}(pt)=\MU_*/(\Sigma).$$  

By Mironov \cite{MI} (Theorem 4.3 and Theorem 4.5) $\MU^*_{\Sigma}(-)$ admits an associate multiplication and all obstructions to commutativity are in $\MU_*/(\Sigma)\otimes \mathbb{F}_2$. Therefore after localization away from 2 
all obstructions vanish and we get a commutative cohomology $$\emph{h}_{\Sigma}^*(-):=\MU^*_{\Sigma}[1/2](-).$$

It is clear that $\emph{h}_{\mathcal{B}}^*(-)$ is complex oriented as the  Atiyah-Hirzebruch
spectral sequence 
$H^*(-,\emph{h}_{\mathcal{B}}^*(pt))\Rightarrow \emph{h}_{\mathcal{B}}^*(-)$
collapses for $\BU(1)\times \BU(1)$. 
\qed

\bigskip

\section{  Related genera on $\MSU_*[1/2]$ and $\MSp_*[1/2]$}

Numerous examples of cobordism theory as a classification tool are considered in the literature. Many examples, such as Spin, Sp, SC and SU cobordism, have a complicated structure, and some have not even been fully calculated yet. However, after localized away from  2,  all the mentioned rings are subrings in $\MU_*[1/2]$ and  polynomial.  

Therefore the restriction of a genus on $\MU_*$ classifying a cohomology $h_{\Sigma}$
defined above gives the corresponding genus on $\MSU[1/2]$ \cite{TO} and $\MSp_*[1/2]$. 

In particular by Buchstaber it is known that 
$\MSp^*[1/2](-)$ can be viewed as a cohomology  theory $\Lambda_{\varkappa^*}$ defined by a projector $\varkappa^*$ 
on $\MU_*[1/2]$ \cite{BU}. Moreover $\MSp_*[1/2]$ forgetful  homomorphisms to $\MU_*[1/2]$
is an isomorphism on  the subring in $\MU_*[1/2]$ generated by the coefficients of the series
$$
p_i(\zeta_1\otimes \zeta_2), ,\,\,\,\zeta_i\to HP^{\infty},\,\,\,i=1,2.
$$

Let $\Lambda^{\perp}$ be the subring in $\MU_*$ generated by the coefficients of the series 
$$
G(x)=c_1(\xi+\bar{\xi})=\sum a_ic_2^i(\xi+\bar{\xi}), \,\,\,\xi\to CP^{\infty}.
$$ 

In \cite{BU-NO} it is proved that 

$$
\Lambda\otimes \Lambda^{\perp}[1/2]=\MU_*[1/2].
$$

Here by using the splitting principle everything can be easily calculated in terms of the universal formal group law. Therefore it is reasonable to give explicit calculation of all above mentioned genera.

Note that in \cite{B-M}  A. Baker and J. Morava constructed another but related projector 
on $\MU[1/2]$ also defining $\MSp[1/2]$, but not fixing its forgetful image in $\MU[1/2]$.

There are other projectors on $\MU[1/2]$ which are not fix the forgetful image maps  but still it is possible to work with in the above sense. For example in \cite{BU-1} Buchstaber constructed the projectors in unitary cobordisms that are related to SU-theory.  In \cite{B5} and \cite{C-P} this result is used for explicit calculations.


\begin{thebibliography}{9}
	
	
\bibitem{BA} N. Baas, {\emph { On bordism theory of manifolds with singularities,}} Math. Scand., {\bf 33}(1973), 279-302.


\bibitem{B-M} A. Baker, J. Morava {\emph {$MSp$ localized away from 2 and odd formal group laws}}, 2014, arXiv:1403.2596v1 [math.AT].





\bibitem{B2} M. Bakuradze, {\emph {Computing the Krichever genus}}, J.  Homotopy Relat. Struct, {\bf 9, 1}(2014), 85-93.

\bibitem{B4} M. Bakuradze, {\emph {Cohomological realization of the Buchstaber
formal group law}}, Uspekhi Mat. Nauk, {\bf 77:5(467)} (2022), 189-190. 

\bibitem{B5} M. Bakuradze, {\emph {Polynomial generators of $\MSU_*[1/2]$ related to classifying maps of certain formal group laws}},  HHA, accepted. 




\bibitem{BU-1} V. M. Buchstaber, {\emph {Projectors in unitary cobordisms that are related to SU-theory}}, Uspekhi Mat. Nauk, 
{\bf 27:6(168)}(1972), 231-232.






\bibitem{BU-U} V. Buchstaber, K. Ustinov, {\emph {Coefficient rings of Buchstaber formal group laws}}, Math. Notices, {\bf 206}:11(2015), 19-60.

\bibitem{B-KH} Bukhshtaber V.M., Kholodov A.N., {\emph {Formal groups, functional equations, and generalized cohomology theories, (English. Russian original)}}, Math. USSR Sbornik {\bf 69:1} (1991), 77–97.

\bibitem{BRUNS} W. Bruns, J. Herzog, {\emph {Cohen-Macaulay rings}}, Cambridge University Press, (1993).


\bibitem{BUSATO} Ph. Busato, {\emph {Realization of Abel’s universal formal group law}}, Math. Z., {\bf 239}(2002), 527–561.


 
\bibitem{C-P} G. Chernykh, T. Panov, {\emph {$SU$-linear operations in complex cobordism and the $c_1$-spherical bordism theory}}, Izvestiya Ross. Akad. Nauk Ser. Mat. {\bf 87} (2023), no. 4.









\bibitem{MI} O. K. Mironov, {\emph {Multiplications in cobordism theories with singularities, and
Steenrod – tom Dieck operations}},  Math. USSR
Izvestija {\bf 13}(1979), 89–106; 


\bibitem{L-C-P} I. Yu. Limonchenko, T. E. Panov, and G. S. Chernykh {\emph {SU-bordism: structure results and geometric representatives}}, Russian Math. Surv., {\bf 74:3}(2019), 461–524.

\bibitem{NOV} S. P. Novikov, {\emph {Homotopy properties of Thom complexes (Russian)}}, Mat. Sb. {\bf 57}(1962), 407–442.



\bibitem{Stong} R. E. Stong, {\emph {Notes On Cobordism Theorey}}, Princeton University Press and University of Tokyo Press, (1968)






\bibitem{TO} B. Totaro, {\emph {Chern numbers for singular varieties and elliptic homology}}, Annals of Mathematics, {\bf 151}(2000), 757–792.

  


 
 


\end{thebibliography}
\end{document}